\documentclass[twocolumn,10pt]{article}

\usepackage[english]{babel}

\usepackage[letterpaper,top=2cm,bottom=2cm,left=3cm,right=3cm,marginparwidth=1.75cm]{geometry}

\usepackage{amsmath}
\usepackage{graphicx}
\usepackage{xurl}
\usepackage[colorlinks=true, allcolors=blue]{hyperref}
\usepackage[numbers]{natbib}
\usepackage{amsmath,amsfonts}
\usepackage{algorithmic}
\usepackage{algorithm}
\usepackage{array}
\usepackage{stmaryrd}
\SetSymbolFont{stmry}{bold}{U}{stmry}{m}{n}
\usepackage[OT1]{fontenc}
\usepackage{amssymb}
\usepackage{amsmath}
\usepackage{ dsfont }
\usepackage{float}
\usepackage{empheq}
\usepackage[caption=false,font=normalsize,labelfont=sf,textfont=sf]{subfig}
\usepackage{textcomp}
\usepackage{url}
\usepackage{verbatim}
\usepackage{graphicx}
\usepackage{xcolor} 

\newcommand{\dint}{\displaystyle\int}

\newcommand{\mat}[1]{\bar{\boldsymbol{#1}}}
\newcommand{\vct}[1]{\boldsymbol{#1}}
\newcommand{\vctu}[1]{\hat{\boldsymbol{#1}}}

\newcommand{\inc}{_{\textrm{inc}}}
\newcommand{\sca}{_{\textrm{sca}}}

\title{Exact Fourier Coupling Relations for the Method of Auxiliary Sources on Elliptic Cylinders}
\author{Prisca {Le Dily}, Arnaud Coatanhay, Marie Kintzinger, and Yves-Marie Scolan}
\date{}

\begin{document}
\maketitle

\begin{abstract}
The Method of Auxiliary Sources has been extensively studied in the case of cylindrical geometries. In this article, the authors provide analysis elements for the elliptical geometry case. Based on a Mathieu functions expansion, the impedance matrix is projected into a reduced diagonal matrix. At the same time, the matrix is studied in the Fourier domain. A checkerboard structure is found and a relation between the different elements of the matrix is proven. Numerical simulations confirm the theoretical results. We additionally identify an exact finite-grid parity rule, an aliasing identity connecting the discrete and continuous Fourier representations, and controlled near-circular Fourier-band localization. These results are mostly independent of a small-$q$ truncation of radial Mathieu functions.\end{abstract}

\section{Introduction}
\begin{figure}
\centering
\includegraphics[width=0.45\textwidth]{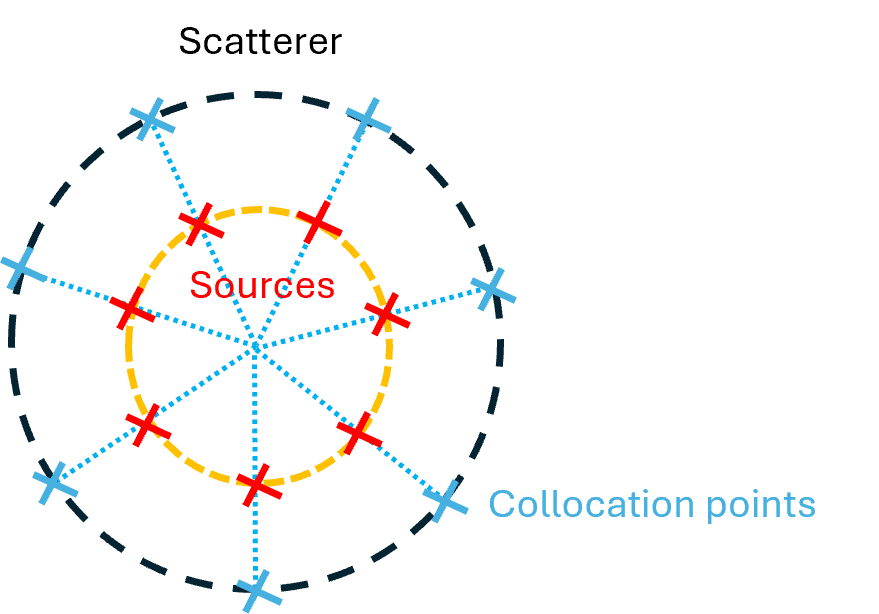} 
\caption{Method of Auxiliary Sources applied to a scatterer, shaped in this case as a cylinder.}\label{fig:MAS}
\end{figure}

The Method of Auxiliary Sources (MAS) is also known as the Method of Fundamental Solutions (MFS), the Discrete Sources Method (DSM) or the desingularized technique. Similarly to the Method of Moments (MoM), it consists of meshing a scatterer and solving boundary conditions at the collocation points, on the scatterer surface. The main difference is that elementary source points are not located at the collocation points, but are displaced inside the scatterer, on a fictitious auxiliary surface (see fig. \ref{fig:MAS}). 
Since singularities occur at the source points, the scattered field needs to be evaluated on those locations for the MoM. On the contrary, for a perfectly conducting scatterer, it is known everywhere for the MAS. This method comes from the calibration procedure of the Boundary Element Method \cite{koopmann_method_1989}. It is widely used in acoustics \cite{koopmann_method_1989,karageorghis_application_2011,doicu_acoustic_2000}, electromagnetism \cite{wei_singular_2019,costa2026aaasquaressolutionhelmholtz,barnett_stability_2008} and hydrodynamics \cite{scolan_aspects_2010}. In the electromagnetism case, a recent and detailed review has been published by \citet{papakanellos_method_2024}, updating a previous review \cite{Kaklamani2002APM}.\newline
The MAS is often used for scatterers with a closed surface \cite{kouroublakis_time-domain_2026,anioutine_about_2003,anastassiu_optimization_2005}, although a few works also involve open surfaces \cite{petoev_method_2015,kou_study_2025,papakanellos_convergence_2003}. The most studied configuration is the canonical case of a circular cylinder \cite{kouroublakis_time-domain_2026,wriedt_generalized_2018,Fikioris2007TAP}. The reason is that it can be easily compared to the exact solution of \citet{rayleigh_xli_1918}. The influence of parameters such as the distance between the source points and the collocation points, or the number of sources, on the error and condition number is a well-known result from \citet{anastassiu_error_2004}. By extending the cylindrical case to a slightly perturbed case, a preconditioning method for ellipses has also been created  by
\citet{anastassiu_efficient_2008}.\newline
As a matter of fact, with the appropriate conformal mapping, any scatterer with a closed surface can be transformed back to the case of a circular cylinder (see the Riemann mapping theorem \cite{nehari2012conformal}). However, conformal maps $\Phi$ do not preserve the constant-coefficient Helmholtz equation: the transformed equation has a spatially varying factor $k_0^2|\Phi^{\prime}|^2$.\newline
As for ellipses, the Joukowski map provides a convenient parametrization of confocal source curves; it is used here to study the Fourier structure of the actual Helmholtz kernel, without replacing the wave problem by a circular Helmholtz problem.\newline
Ellipses are a canonical case as well, as an exact solution exists \cite{van_den_berg_diffraction_1973}, based on a decomposition into Mathieu functions.\newline
Similarly, the elements of the MAS impedance matrix can be decomposed into sums of Mathieu functions. Following this observation, a Singular Value Decomposition (SVD) was proposed to improve the conditioning of the MAS \cite{antunes_improving_2024}.\newline
The purpose of this article is to highlight a few analysis elements on the Method of Auxiliary Sources for ellipses, based on the decomposition of the impedance matrix into Mathieu functions. The sampled matrix admits an approximately diagonal reduced modal compression, with explicit continuous Mathieu multipliers. This analysis does not require a small eccentricity.\newline
Furthermore, the cylindrical case exhibited a diagonal impe\-dance matrix in the Fourier basis. This article presents an analysis of the impedance matrix behavior in the Fourier basis, for the elliptical case.\newline
After a brief reminder of the MAS procedure, the impedance matrix is decomposed and reduced to a smaller diagonal matrix. Then, the impedance matrix is also expressed in the Fourier basis, which reveals a checkerboard pattern. A relation between the different matrix elements is established. We further separate exact continuous identities from discrete sampling effects and establish a near-circular band-localization estimate. Finally, numerical results confirm our observations.

\section{Methodology}
In this communication, a unit vector is denoted with a hat $\vctu{u}$, a vector is bold $\vct{u}$, and their norm is $u$. The time convention $e^{-j\omega t}$ is applied. The ellipse surface is perfectly-conducting. The incident field is a plane wave $\psi\inc(\vct{r})=e^{j\vct{k}\inc\cdot\vct{r}}$, with $\vct{k}\inc$ the wave vector and $\vct{r}$ the position. The wavenumber is named $k_0$.
\subsection{The Method of Auxiliary Sources (MAS)}
The scattered field $\psi\sca$ is considered at the observation point $\vct{r}\sca$. On the ellipse surface, collocation points are denoted $\vct{r}_k$ and the source points are $\vct{r}_{s,k}$ with $k\in\llbracket 0;P-1\rrbracket$. At any spatial coordinate, the scattered field is written with the MAS as
\begin{align}\label{2DMAS}
    \psi\sca(\vct{r}\sca)=\sum_{i=0}^{P-1} \alpha_i g\left(\vct{r}\sca,\vct{r}_{s,i}\right)
\end{align}
with $g\left(\vct{r}\sca,\vct{r}_{s,i}\right)$ the Green function and $P$ the number of sources. For a 2D-geometry, the Green function is \begin{equation}
    g\left(\vct{r}\sca,\vct{r}_{s,i}\right)=\frac{j}{4} H_0^{(1)}\left(k_0||\vct{r}\sca-\vct{r}_{s,i} ||\right)
\end{equation}
Coefficients $\alpha_i$ are determined by applying the boundary conditions at the scatterer surface.\newline

For the scalar Dirichlet problem considered here, the boundary condition is 
\begin{equation}\label{BC}
    \forall k,\, \psi\inc(\vct{r}_k)+\psi\sca(\vct{r}_k) = 0
\end{equation}
with
\begin{align}\label{2Dres}
    \psi\sca(\vct{r}_k)=\frac{j}{4}\sum_{i=0}^{P-1} \alpha_i  H_0^{(1)}\left(k_0||\vct{r}_k-\vct{r}_{s,i} ||\right)
\end{align}
This leads us to the resolution of a linear system of the type $\mat{Z}\vct{X}=\vct{B}$, where $\vct{X}$ is the vector $(\alpha_i)_i$ and $\vct{B}$ the incident field vector $(-\psi_\textrm{inc,k})_k$.
$\mat{Z}$ is the impedance matrix of dimensions $P\times P$, written as 
\begin{align}\label{ped}\centering
&\mat{Z}=\frac{j}{4}\times\nonumber\\
&\resizebox{.88\hsize}{!}{$\begin{pmatrix}H_0^{(1)}\left(k_0||\vct{r}_0-\vct{r}_{s,0} ||\right) &\dots &H_0^{(1)}\left(k_0||\vct{r}_0-\vct{r}_{s,P-1} ||\right)\\
    \vdots&\ddots&\vdots\\
    H_0^{(1)}\left(k_0||\vct{r}_{P-1}-\vct{r}_{s,0} ||\right) &\dots &H_0^{(1)}\left(k_0||\vct{r}_{P-1}-\vct{r}_{s,P-1} ||\right) \end{pmatrix}$}
\end{align}
Elliptical coordinates on the ellipse surface are $\vct{r}_k = (x_k,y_k)$, with
\begin{align}\label{ellipco}
    x_k &= h\cosh\mu_s\cos\nu_k\nonumber\\
    y_k &= h\sinh\mu_s \sin\nu_k
\end{align}
for $P$ collocation points, positioned at uniformly distributed angles $\nu_k=\frac{2\pi k}{P}$, where $k\in\llbracket 0;P-1\rrbracket$.\newline
The convergence of the MAS is not guaranteed if singularities are not enclosed by the fictitious auxiliary surface \cite{kyurkchan_singularities_1999,Tsitsas2018TAP}. For elliptical scatterers, the focal segment is a natural analytic-continuation feature; a confocal auxiliary ellipse enclosing it provides a geometrically convenient source curve, without by itself guaranteeing convergence. This is the reason why sources are located on a confocal ellipse (i.e., with the same focal half-distance $h$ as the elliptical scatterer), with $\vct{r}_{s,k} = (x_{s,k},y_{s,k})$, and
\begin{align}\label{ellipcos}
    x_{s,k} &= h\cosh\mu_\delta\cos\nu_{s,k}\nonumber\\
    y_{s,k} &= h\sinh\mu_\delta \sin\nu_{s,k}
\end{align}
where $0<\mu_\delta<\mu_s$.  \newline
The impedance matrix can be decomposed into a sum of Mathieu functions. As a matter of fact, according to \cite{betcke_spectral_2014,morse1953methods}, 
\begin{align}\label{hankmat}
&\frac{j}{4} H_0^{(1)}\big(k_0\|\vct{r}_1 - \vct{r}_2\|\big)=\frac{j}{2}\times
\nonumber\\
&\Bigg[\sum_{m=0}^{\infty}
\mathrm{ce}_m(\nu, q)\,\mathrm{ce}_m(\nu', q)\,
\mathrm{Mc}_m^{(1)}(\mu_<, q)\,
\mathrm{Mc}_m^{(3)}(\mu_>, q)
\nonumber\\
&
+
\sum_{m=1}^{\infty}
\mathrm{se}_m(\nu, q)\,\mathrm{se}_m(\nu', q)\,
\mathrm{Ms}_m^{(1)}(\mu_<, q)\,
\mathrm{Ms}_m^{(3)}(\mu_>, q)
\Bigg]
\end{align}

The Mathieu parameter is $q=\frac{k_0^2h^2}{4}$. $(\mu,\nu)$ and $(\mu',\nu')$ are the respective elliptical coordinates of $\vct{r}_1$ and $\vct{r}_2$. $\mu_<,\mu_>$ correspond respectively to the smallest and highest coordinate $\mu$ of points $\vct{r}_1$ and $\vct{r}_2$.
The decomposition \eqref{hankmat} consists of two sums: cosine-even and sine-odd Mathieu families. The radial factors and overall prefactor in \eqref{hankmat} are understood in the normalization used here; changing the normalization of angular Mathieu functions requires a corresponding modification of the radial coefficients.\newline
Let $\vct{r}_1$ be linked to the collocation points $(\mu_k,\nu_k)$ and $\vct{r}_2$ to the source points $(\mu_\delta,\nu_{s})$. The impedance matrix \eqref{ped} can be split into two matrices. The first matrix $\mat{Z}_1$ is constituted of even elements,
\begin{align}\label{Z_1}
 &\mat{Z}_1 = \frac{j}{2}\times\nonumber\\
 &\resizebox{.89\hsize}{!}{$\left( \sum_{m=0}^{\infty}
\mathrm{ce}_m(\nu_k, q)\,\mathrm{ce}_m(\nu_s, q)\,
\mathrm{Mc}_m^{(1)}(\mu_\delta, q)\,
\mathrm{Mc}_m^{(3)}(\mu_s, q) \right)_{k,s}$}
\end{align}

and the second $\mat{Z}_2$ consists of the odd Mathieu modes
\begin{align}\label{Z_2}
& \mat{Z}_2= \frac{j}{2}\times\nonumber\\
&\resizebox{.89\hsize}{!}{$\left( \sum_{m=1}^{\infty}
\mathrm{se}_m(\nu_k, q)\,\mathrm{se}_m(\nu_s, q)\,
\mathrm{Ms}_m^{(1)}(\mu_\delta, q)\,
\mathrm{Ms}_m^{(3)}(\mu_s, q) \right)_{k,s}$}
\end{align}
Both matrices are complex symmetric (non-Hermitian). To simplify the calculations, we assume that collocation and source points are angularly aligned for $j\in\llbracket0;P-1\rrbracket$.
\subsection{Analysis of the impedance matrix}

Mathieu functions are orthogonal in $\mathds{L}^2$. The propagation kernel $g$ between two confocal ellipses is given by equation \eqref{hankmat}, with $g(\nu_,\nu')$ denoting $g(\vct{r}(\mu_s,\nu),\vct{r}(\mu_\delta,\nu'))$. In the even case, we can notice that
\begin{align}\label{integr}
&   \dint_0^{2\pi}\dint_0^{2\pi}
g(\nu,\nu')\, \mathrm{ce}_n(\nu, q)\,\mathrm{ce}_l(\nu', q) d\nu d\nu'\nonumber\\
&= \left\{
    \begin{array}{ll}
        0 & \mbox{if } n\neq l \\
        \frac{j\pi^2}{2} \mathrm{Mc}_n^{(1)}(\mu_\delta, q)\,
\mathrm{Mc}_n^{(3)}(\mu_s, q) & \mbox{otherwise.}
    \end{array}
\right.
\end{align}
The integral $\dint_0^{2\pi}\mathrm{ce}_n(\nu, q)\,\mathrm{ce}(\nu, q)d\nu=\delta_{n,m}\pi$ can be approached by 
$
    \frac{2\pi}{P}\sum_{i=0}^{P-1} \mathrm{ce}_n(\nu_i, q)\,\mathrm{ce}_m(\nu_i, q)\approx \delta_{n,m}\pi
$. \newline

Equation \eqref{integr} can be approached by 
\begin{equation}
    \left(\frac{2\pi}{P}\right)^2  \sum_{k=0}^{P-1}\sum_{s=0}^{P-1}
\mat{Z}_{k,s} \mathrm{ce}_n(\nu_k, q)\,\mathrm{ce}_l(\nu_s, q) 
\end{equation}
This can be formulated as a matrix product :
\begin{equation}
 \left(\frac{2\pi}{P}\right)^2  \left(\mat{P}_1^T \mat{Z}_1 \mat{P}_1\right)_{n,l}
\end{equation}
with
\begin{align}\label{passMat}
    \mat{P}_1=\begin{pmatrix}
    \mathrm{ce}_0(\nu_0)&\dots&\mathrm{ce}_{M-1}(\nu_0)\\
    \vdots&\ddots&\vdots\\
        \mathrm{ce}_0(\nu_{P-1})&\dots&\mathrm{ce}_{M-1}(\nu_{P-1})\\
    \end{pmatrix}
\end{align}
The matrix columns can be normalized, through a multiplication by $\sqrt{\frac{2\pi}{P}}\times\frac{1}{\sqrt{\pi}}=\sqrt{\frac{2}{P}}$. The resulting $P\times M$ matrix can have orthonormal columns in principle; for finite $P$ their orthogonality is only approximate because the trapezoidal rule samples Mathieu functions. The reduction is a modal compression, not an exact similarity diagonalization of the full collocation matrix.\newline
To simplify the calculations, $\mat{Z}_1$ is also normalized by $\frac{2}{P}$.
The compressed, rescaled matrix $\mathbf{\tilde{P}}_1^T \mathbf{\tilde{Z}}_1 \mathbf{\tilde{P}}_1$ is approximately diagonal for resolved modes. Its expected modal diagonal entries (with $\widetilde{Z}_1=\frac{2}{P}\mat{Z}_1$) are
\begin{equation}\label{vpcm}
 \boxed{ \lambda_n =  \frac{j}{2} \mathrm{Mc}_n^{(1)}(\mu_\delta, q)\,
\mathrm{Mc}_n^{(3)}(\mu_s, q)}, \qquad n\geq 0.
\end{equation}
Let us consider $\mat{Z}_2$ in the same way. The change-of-basis matrix is this time
\begin{align}\label{passMat2}
    \mat{P}_2=\begin{pmatrix}
    \mathrm{se}_1(\nu_0)&\dots&\mathrm{se}_{M}(\nu_0)\\
    \vdots&\ddots&\vdots\\
        \mathrm{se}_1(\nu_{P-1})&\dots&\mathrm{se}_{M}(\nu_{P-1})\\
    \end{pmatrix}
\end{align}
Columns are normalized in the same way. Similarly, the expected modal diagonal values are
\begin{equation}\label{vpsm}
 \boxed{ \lambda_n =  \frac{j}{2} \mathrm{Ms}_n^{(1)}(\mu_\delta, q)\,
\mathrm{Ms}_n^{(3)}(\mu_s, q)}, \qquad n\geq 1.
\end{equation}
\subsection{Expression in the Fourier basis}
The angular Mathieu functions can be expressed as (see 20.2.3 and 20.2.4 of \citet{abramowitz_handbook_1948}) :
\begin{equation}\label{ce_n}
    \mathrm{ce}_{2m+p}(z,q) = \sum_{n=0}^{\infty} \mathcal{A}_{2n+p}^{\,2m+p}(q) \cos\big((2n+p)z\big), \quad p=0,1,
\end{equation}
and
\begin{equation}\label{se_n}
    \mathrm{se}_{2m+p}(z,q) = \sum_{n=0}^{\infty} \mathcal{B}_{2n+p}^{\,2m+p}(q) \sin\big((2n+p)z\big), \quad p=0,1.
    \end{equation}
In equation \eqref{se_n}, the first coefficient of the series is null.
The Fourier change-of-basis matrix is defined by $\mat{F}_{k,s}=\frac{1}{\sqrt{P}}e^{js\nu_k}$, while $\mat{F}^\dag_{l,s}=\frac{1}{\sqrt{P}}e^{-il\nu_s}$.\newline
With $L^m_{\mu_s,\mu_\delta}= \mathrm{Mc}_m^{(1)}(\mu_\delta, q)\,
\mathrm{Mc}_m^{(3)}(\mu_s, q) $, this leads to
\begin{align}\label{prodmat}
  &  \left(\mat{Z}_{1}\mat{F}\right)_{l,s}\nonumber\\
&=\frac{j}{2\sqrt{P}}\sum_{k=0}^{P-1} \sum_{m=0}^{\infty}L^m_{\mu_s,\mu_\delta}
\mathrm{ce}_m(\nu_l, q)\,\mathrm{ce}_m(\nu_k, q)\,e^{js\nu_k}\,
\nonumber\\
&=\frac{j}{2\sqrt{P}} \sum_{m=0}^{\infty}L^m_{\mu_s,\mu_\delta}
\mathrm{ce}_m(\nu_l, q)\,
\sum_{k=0}^{P-1}\mathrm{ce}_m(\nu_k, q)\,e^{js\nu_k}
\end{align}
For the complex Fourier representation, define $\mathcal C_r^{m}=\mathcal A_{|r|}^{m}$ when $r\ne0$, and $\mathcal C_0^{m}=2\mathcal A_0^{m}$. The doubled zero harmonic is required because $\cos(0\nu)=1$. According to \eqref{ce_n}, $\mathrm{ce}_n(\nu,q)=\frac{1}{2}\sum_{r=-\infty}^{\infty} \mathcal C_r^{n}(q) e^{jr\nu}$.
Then,
\begin{align}\label{Fouriercem}
    \sum_{k=0}^{P-1}\mathrm{ce}_m(\nu_k, q)\,e^{js\nu_k} =  \frac{1}{2}\sum_{k=0}^{P-1}\sum_{r=-\infty}^\infty \mathcal{C}_r^{m}(q) e^{jr\nu_k}\,e^{js\nu_k}
\end{align}
Knowing that $\nu_k=\frac{2\pi k}{P}$,
 \begin{align}\label{sumexp}
      \sum_{k=0}^{P-1} e^{j(r+s)\nu_k}&= \left\{
    \begin{array}{ll}
        \frac{1-1}{1-e^{j(r+s)\frac{2\pi }{P}}} & \mbox{if } r\neq -s+nP \\
        P & \mbox{otherwise.}
    \end{array}
\right.
 \end{align}
 It means that \begin{align}
    \sum_{k=0}^{P-1}\mathrm{ce}_m(\nu_k, q)\,e^{js\nu_k} =   \frac{P}{2}\sum_{n\in\mathbb{Z}} \mathcal{C}_{-s+nP}^{m}(q) 
\end{align}
In numerical terms, the calculation of $\mathrm{ce}_n=\frac{1}{2}\sum_{r=-R}^{R}\mathcal{C}_r^{n} e^{jr\nu}$ has to be a finite sum. For a Fourier truncation $|r|\le R$, choosing $P>2R$ suppresses aliasing of that truncated series. Since an exact Mathieu function is not band-limited, the following identity holds only in this alias-free truncation (or asymptotically for fixed resolved $s$ as $P$ increases):
\begin{align}\label{ressom}
    \sum_{k=0}^{P-1}\mathrm{ce}_m(\nu_k, q)\,e^{js\nu_k} \simeq   \frac{P}{2} \mathcal{C}_{s}^{m}(q) 
\end{align}
since $\mathcal C_{-s}^{m}=\mathcal C_s^{m}$. As a consequence,
\begin{align}
    \left(\mat{Z}_{1}\mat{F}\right)_{l,s}
&=\frac{j\sqrt{P}}{4} \sum_{m=0}^{\infty} 
L^m_{\mu_s,\mu_\delta}\mathrm{ce}_m(\nu_l, q)\,
\mathcal{C}_{s}^{m}(q)\,
\end{align}
Then,
\begin{align}
&\left(\mat{F}^\dag\mat{Z}_{1}\mat{F}\right)_{k,s}\nonumber\\
&=\frac{j}{4} \sum_{m=0}^{\infty}
L^m_{\mu_s,\mu_\delta}\mathcal{C}_{s}^{m}(q)\,\sum_{l=0}^{P-1}\mathrm{ce}_m(\nu_l, q)\,e^{-ik\nu_l}
\end{align}
Similarly to equation \eqref{ressom},
 \begin{align}
    \sum_{l=0}^{P-1}\mathrm{ce}_m(\nu_l, q)\,e^{-ik\nu_l} \simeq   \frac{P}{2}\mathcal{C}_{k}^{m}(q)
\end{align}
Finally,
\begin{align}\label{sommedam}\boxed{
\left(\mat{F}^\dag\mat{Z}_{1}\mat{F}\right)_{k,s}\simeq\frac{jP}{8} \sum_{m=0}^{\infty}
L^m_{\mu_s,\mu_\delta}\mathcal{C}_{s}^{m}(q)\,\mathcal{C}_{k}^{m}(q)}
\end{align}
This expression depends on the product $\mathcal{C}_{s}^{m}(q)\,\mathcal{C}_{k}^{m}(q)$. We know that if $m$ is even, only even indices $s$ and $k$ are nonzero. The same applies for $m$ odd. As a consequence, the product $\mathcal{C}_{s}^{m}(q)\,\mathcal{C}_{k}^{m}(q)$ is non zero only if $s$ and $k$ have the same parity. It means that $\hat{\mathbf{Z}}_1=\mat{F}^\dag\mat{Z}_{1}\mat{F}$ has a checkerboard structure, as every second element is null.\newline

Let us now consider the second matrix $\mat{Z}_2$ from equation \eqref{Z_2}. According to \eqref{se_n}, 
\begin{equation}
    \mathrm{se}_n(\nu,q)=\frac{1}{2j}\sum_{r\in\mathbb Z}\textrm{sgn}(r)\, \mathcal{B}_{|r|}^{n}(q) \,e^{jr\nu}
\end{equation} 
The calculations are very similar to the even case. Let us denote $K^m_{\mu_s,\mu_\delta}= \mathrm{Ms}_m^{(1)}(\mu_\delta, q)\,
\mathrm{Ms}_m^{(3)}(\mu_s, q) $. We also have
\begin{align}
    \sum_{k=0}^{P-1}\mathrm{se}_m(\nu_k, q)\,e^{js\nu_k} \simeq \frac{jP}{2} \textrm{sgn}(s)\,\mathcal{B}_{|s|}^{m}(q) 
\end{align}
while
 \begin{align}
    \sum_{l=0}^{P-1}\mathrm{se}_m(\nu_l, q)\,e^{-ik\nu_l} \simeq -\frac{jP}{2}\textrm{sgn}(k)\,\mathcal{B}_{|k|}^{m}(q)
\end{align}

Finally,
\begin{empheq}[box=\fbox]{align}\label{sommedam2}
 &\left(\mat{F}^\dagger \mat{Z}_2F\right)_{k,s}
 \simeq \frac{jP}{8}\operatorname{sgn}(s)\operatorname{sgn}(k)\nonumber\\
 &\hspace{1em}\times\sum_{m=0}^{\infty}
 K^m_{\mu_s,\mu_\delta}\mathcal B_{|s|}^{m}(q)
 \mathcal B_{|k|}^{m}(q).
\end{empheq}
Mathieu parity gives the same checkerboard rule for $ \hat{\textbf{ Z}}_2=\mat{F}^\dagger \mat{Z}_2 F$. Hence the full matrix $\hat{\textbf{ Z}}=\hat{\textbf{ Z}}_1+\hat{\textbf{ Z}}_2$ obeys the same selection rule.\par
The factor $1/j$ in the sine expansion is essential: at $q=0$, the cosine and sine sectors must add to the circular Fourier diagonal, rather than cancel. The approximate signs above neglect Fourier aliases for finite $P$.\par
\subsection{Relations between elements in the Fourier basis}
The Mathieu equation, whose solutions are $\mathrm{ce}_n$ and $\mathrm{se}_n$, is
\begin{equation} \frac{d^2y}{d\nu^2} + (a - 2q\cos 2\nu)\,y = 0. 
\end{equation}
Let us introduce the Mathieu angular operator \begin{equation} L_q=-\frac{d^2}{d\nu^2}+2q\cos 2\nu. \label{20} \end{equation} 
The functions $\mathrm{ce}_n$ and $\mathrm{se}_n$ are eigenfunctions of $L_q$, with their respective characteristic values $a_n(q)$ and $b_n(q)$. \newline
An angular integral operator associated with the kernel from equation \eqref{hankmat} can be defined, at the continuous level, as  \begin{equation} (Tf)(\nu_s)=\int_0^{2\pi} g(\nu_s,\nu_\delta)\,f(\nu_\delta)\,d\nu_\delta, \label{21} \end{equation}
Its Green expansion in Mathieu functions shows that the same families $\mathrm{ce}_n,\mathrm{se}_n$ diagonalize $T$. With $\int_0^{2\pi}\mathrm{ce}_n^2 d\nu=\int_0^{2\pi}\mathrm{se}_n^2 d\nu=\pi$, the eigenvalues of the integral operator $T$ are $\pi$ times the rescaled modal entries in \eqref{vpcm}--\eqref{vpsm}. Using the operator $L_q$ for $\nu=\nu_s$, one has the intertwining relation $L_qT \mathrm{ce}_n=TL_q\mathrm{ce}_n$. As angular Mathieu functions are a complete set for $\pi$- or $2 \pi$-periodic functions, one has
\begin{equation} 
L_qT=TL_q . \label{22} 
\end{equation} 
for functions in the corresponding periodic domain.\newline
Writing $\widehat{T}_{k,s}=\langle e_k,T e_s\rangle$, $e_k(\nu)=e^{jk\nu}$ and the normalized Fourier inner product \begin{equation}
    \langle f,g\rangle=(2\pi)^{-1}\dint_0^{2\pi}f^*(\nu)g(\nu)\,d\nu,
\end{equation}
with $f^*(\nu)$ the complex conjugate of $f(\nu)$, one has
\begin{equation}
L_q e_k = k^2 e_k + q e_{k+2} + q e_{k-2}. \label{23}
\end{equation}
Taking the $(k,s)$ Fourier matrix element of the continuous commutation relation \eqref{22} therefore yields
\begin{equation}\label{resid}
(k^2-s^2)\,\widehat{T}_{k,s} + q\bigl( \widehat{T}_{k-2,s} + \widehat{T}_{k+2,s} - \widehat{T}_{k,s-2} - \widehat{T}_{k,s+2} \bigr) =0.
\end{equation}
For uniformly sampled angles $\nu_j=2\pi j/P$, the unitary DFT (Discrete Fourier Transform) of the MAS matrix, $\widehat{\mathbf{Z}}^{(P)}=\mat{F}^{\dag}\mat{Z}\mat{F}$, is the trapezoidal discretization of the Fourier representation of $T$: for resolved modes, $\widehat{T}_{k,s}\simeq(2\pi/P)\widehat{\mathbf{Z}}^{(P)}_{k,s}$. Thus, \eqref{resid} is exact at the continuous level, while the corresponding discrete residual converges to zero for each fixed, resolved pair $(k,s)$ as $P\to\infty$, provided the kernel is sampled consistently and the source and boundary ellipses remain separated. This statement does not assert uniform convergence over all Fourier indices up to Nyquist.


\subsection{Exact discrete parity and aliasing}
There is a useful distinction between the exact parity rule and the
five-point relation. Since the two confocal curves are centrally symmetric,
$g(\nu+\pi,\nu'+\pi)=g(\nu,\nu')$. Hence
\begin{equation}\label{parityop}
 \widehat T_{k,s}=0\qquad\text{if }k-s\text{ is odd}.
\end{equation}

For an even number of angular values $P$, the half-turn maps the uniform grid into itself. Therefore
\eqref{parityop} also holds \emph{exactly} for
$\mat{F}^\dagger \mat{Z}_P\mat{F}$, independently of the angular quadrature error.
For odd $P$, this exact finite-grid argument does not apply.

For the unitary discrete Fourier transform, introduce
$\tilde{\textbf{ Z}}^{(P)}=(2\pi/P)\mat{F}^\dagger \mat{Z}_P\mat{F}$.
Absolute convergence of the double Fourier series of the smooth kernel
gives the exact aliasing identity
\begin{equation}\label{exactalias}
 \tilde{\textbf{ Z}}^{(P)}_{k,s}=
 \sum_{a,b\in\mathbb Z}\widehat T_{k+aP,s+bP}.
\end{equation}
Thus finite-grid formulas obtained from truncated Mathieu--Fourier
series must retain the missing aliases, rather than being treated as
exact identities for the full sampled matrix.

\subsection{Localization of Fourier couplings near a circle}
The source-curve location and the shape deformation can be parametrized
independently using the Joukowski map. With $a$ the major semi-axis and $b$ the minor semi-axis,
\begin{equation}\label{JoukowskiEabe}
 \Phi_\varepsilon(\zeta)=c_s(\zeta+\varepsilon\zeta^{-1}),\quad
 c_s=\frac{a+b}{2},\quad
 \varepsilon=\frac{a-b}{a+b},
\end{equation}
where the scatterer corresponds to $|\zeta|=1$ and the confocal
source curve to $|\zeta|=\rho$ with
$\sqrt\varepsilon<\rho<1$. The precise relation to elliptic
coordinates is $\rho=\exp(\mu_\delta-\mu_s)$, and
$q=(k_0c_s)^2\varepsilon$. In particular, $\varepsilon$ is a
geometric deformation parameter, whereas $q$ also involves frequency.
Set $u=\nu-\nu'$ and $v=\nu+\nu'$. A direct calculation gives
\begin{align}\label{JoukdistanceEabe}
 \frac{|\Phi_\varepsilon(e^{j\nu})-
 \Phi_\varepsilon(\rho e^{j\nu'})|^2}{c_s^2}
 ={}& A_0(u)+2\varepsilon B_1(u)\cos v\nonumber\\
 &+\varepsilon^2 A_2(u),\\
 A_0(u)={}&1+\rho^2-2\rho\cos u,\nonumber\\
 B_1(u)={}&2\cos u-(\rho+\rho^{-1}),\nonumber\\
 A_2(u)={}&1+\rho^{-2}-2\rho^{-1}\cos u.\nonumber
\end{align}

At fixed electrical size $k_0c_s$ and fixed source ratio $\rho\in(0,1)$,
the Helmholtz Green kernel is analytic in $\varepsilon$ in a
neighborhood of $0$, because the two curves remain disjoint.\newline
Let us denote the Green kernel as
$K_\varepsilon(u,v)
=
F\bigl(D_\varepsilon(u,v)\bigr),
$ with the function
$F(x)
=
\frac{j}{4}
H_0^{(1)}\bigl(k_0c_s\sqrt{x}\bigr)
$ and  $D_\varepsilon(u,v) = \frac{|\Phi_\varepsilon(e^{j\nu})-
 \Phi_\varepsilon(\rho e^{j\nu'})|^2}{c_s^2}$. By denoting
$\delta
=
2\varepsilon B_1(u)\cos v
+\varepsilon^2 A_2(u)$, we can write
\begin{equation}
F(A_0+\delta)
=
\sum_{p=0}^{\infty}
\frac{F^{(p)}(A_0)}{p!}\delta^p.
\end{equation}

Developing $\delta^p$ provides terms of the form \begin{equation}
    C_{a,b}(u)\,\varepsilon^{a+2b}\cos^a v,
 a,b\in\mathbb N
\end{equation}
Terms of order $n$ for $\varepsilon$
must then verify $
n=a+2b$. Moreover,
\begin{equation}
\cos^a v
=
\frac{1}{2^a}
\sum_{i=0}^{a}
\binom{a}{i}e^{j(a-2i)v}.
\end{equation}
With $\ell=a-2i$, we have $|\ell|\leq a\leq n$ and $\ell$ of the same parity as $n$. Writing $K_\varepsilon(u,v)=\sum_{n=0}^{\infty}
\varepsilon^n K_n(u,v)$, each term $K_n$ admits a decomposition of the form
\begin{equation}
K_n(u,v)
=
\sum_{\substack{|\ell|\leq n\\
                 \ell\equiv n\ ({\rm mod}\ 2)}}
k_{n,\ell}(u)e^{j\ell v}.
\end{equation}
Finally, since $e^{jru}e^{j\ell v}=e^{j(r+\ell)\nu}
e^{j(\ell-r)\nu'}$, the elements $(k,s)$ of the Fourier matrix must satisfy $k=r+\ell$ and $s=r-\ell$, and we have
$k-s=2\ell$.
For each couple $(k,s)$, it implies that
\begin{equation}\label{bandhierarchyEabe}
\boxed{
\widehat T_{k,s}(\varepsilon)
=
O\left(\varepsilon^{|k-s|/2}\right).
}
\end{equation}
for each fixed pair $(k,s)$, while opposite parities vanish
exactly. Moreover, if $\mathcal B_R$ deletes all Fourier
bands $|k-s|>2R$, analyticity in the Hilbert--Schmidt
norm yields
\begin{equation}\label{bandtailEabe}
 \|T_\varepsilon-\mathcal B_RT_\varepsilon\|_{\rm HS}
 \le C_R|\varepsilon|^{R+1}
 \quad (|\varepsilon|\text{ sufficiently small}).
\end{equation}
Here $R$ is fixed and the constants may depend on $k_0c_s$
and $\rho$. This is a controlled near-circular band approximation,
\emph{not} a universal bandedness claim for eccentric scatterers.
In particular, numerical comparisons performed at fixed
$\mu_\delta/\mu_s$ need not have fixed $\rho$ and cannot by
themselves test the asymptotic order in \eqref{bandhierarchyEabe}.

\subsection{Confocal separation and Fourier decay}
Let $\Delta\mu=\mu_s-\mu_\delta>0$ and $\Sigma\mu=\mu_s+\mu_\delta$. Based on equations \eqref{ellipco} and \eqref{ellipcos}, the squared distance
between the two confocal curves also admits the exact factorization
\begin{align}\label{ellipdistEabe}
 \frac{|\vct{r}(\mu_s,\nu)-\vct{r}(\mu_\delta,\nu')|^2}{h^2}
 ={}&[\cosh\Delta\mu-\cos(\nu-\nu')]\nonumber\\
 &\times[\cosh\Sigma\mu-\cos(\nu+\nu')].
\end{align}
The nearest complex-angular singularities are therefore separated
from the real torus by $\Delta\mu$ in the difference variable.
For any $\alpha,\beta>0$ with $\alpha+\beta<\Delta\mu$,
contour displacement in the two angular variables gives
\begin{equation}\label{expboundEabe}
 |\widehat T_{k,s}|
 \le C_{\alpha,\beta}
 e^{-\alpha|k|-\beta|s|}.
\end{equation}
Combining this bound with \eqref{exactalias} proves exponential
convergence of the sampled Fourier coefficients to the continuous ones for each fixed mode pair.

\subsection{Conditioning of a truncated Mathieu compression}
For a fixed set of $M$ normalized real angular Mathieu functions
$\varphi_n$ satisfying $\int_0^{2\pi}|\varphi_n|^2=\pi$,
define
\begin{equation}\label{QPEabe}
 (\mat{Q}_P)_{j,n}=\sqrt{\frac{2}{P}}\,\varphi_n(\nu_j),
 \qquad \mat{G}_P=\mat{Q}_P^T\mat{Q}_P.
\end{equation}
The finite modal approximation to the \emph{weighted} collocation
matrix is $\mat{A}_{P,M}=\mat{Q}_P\mat{\Lambda}_M\mat{Q}_P^T$, with $\mat{\Lambda}_M$
the corresponding continuous modal multipliers. If $\mat{G}_P$
is nonsingular and all retained multipliers are nonzero, the
ratio of the largest to the smallest \emph{nonzero} singular
value satisfies
\begin{equation}\label{gramboundEabe}
 \frac{\kappa_2(\mat{\Lambda}_M)}{\kappa_2(\mat{G}_P)}
 \le\kappa_2^+(\mat{A}_{P,M})
 \le\kappa_2(\mat{G}_P)\kappa_2(\mat{\Lambda}_M).
\end{equation}
Indeed, writing $\mat{Q}_P=U_PR_P$ as a thin QR factorization gives
$\mat{A}_{P,M}=U_P(R_P\mat{\Lambda}_MR_P^T)U_P^T$ and
$\kappa_2(R_P)^2=\kappa_2(\mat{G}_P)$.
This estimate concerns the \emph{rank-$M$ truncated model}; it
must not be identified without an additional truncation-error
estimate with the condition number of the full collocation matrix.


\section{Numerical results}
\begin{figure}
\centering
\includegraphics[width=\linewidth]{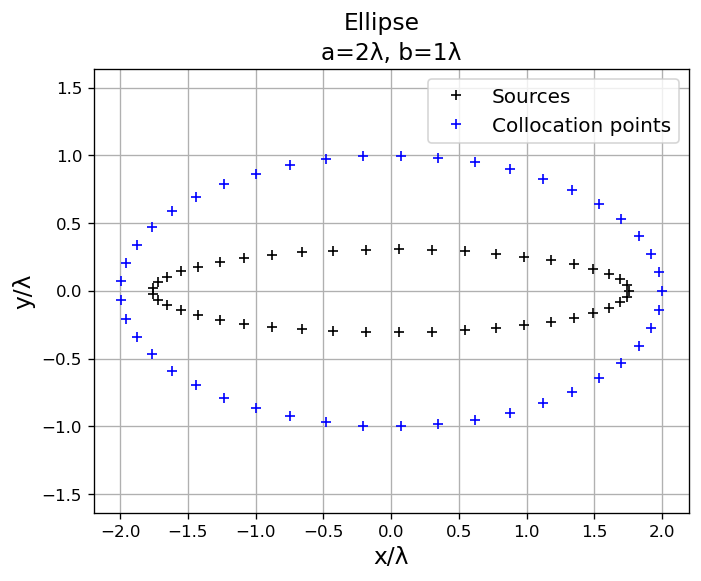} 
\caption{Meshing of the ellipse considered in the numerical simulations.}\label{fig:ellipse}
\end{figure}
\begin{figure}[!t]
\centering
\includegraphics[width=\linewidth]{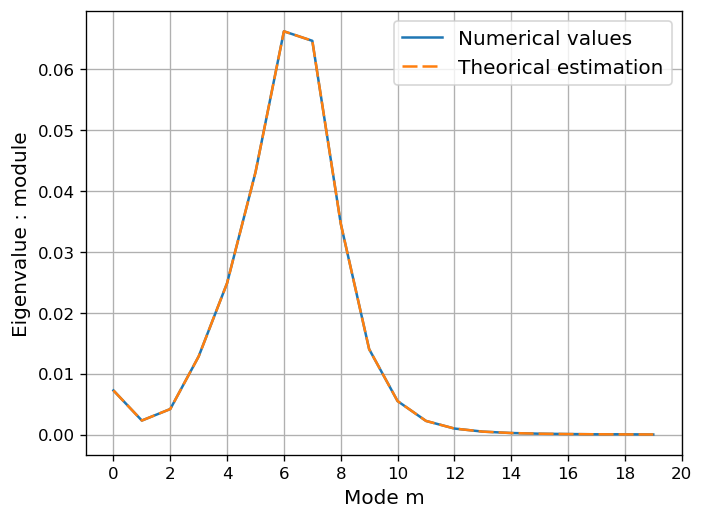}
\caption{Comparison between theoretical and numerical eigenvalues depending on the mode $m$. Even case $\mat{Z}_1$. $\lambda_0=1$ m, $a = 2$ m, $b = 1$ m, $q =29.61$.}
\label{fig:vpcomp}
\end{figure}
\begin{figure}[!t]
\includegraphics[width=\linewidth]{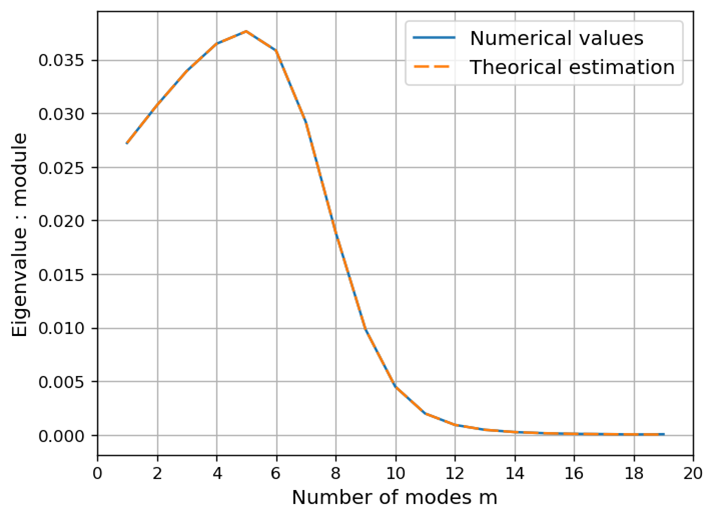}
\caption{Comparison between theoretical and numerical eigenvalues depending on the mode $m$. Odd case $\mat{Z}_2$. $\lambda_0=1$ m, $a = 2$ m, $b = 1$ m, $q =29.61$.}
\label{fig:vpcomp2}
\end{figure}

The previous results can be assessed numerically. In this section, the wavelength is $\lambda_0=1$ m. The major semi-axis of the ellipse is $a=2\lambda_0$, and the minor semi-axis is $b=\lambda_0$. The corresponding Mathieu parameter is $q=\frac{k_0^2\,(a^2-b^2)}{4}=29.61$.\newline
\begin{figure}[!t]
\centering
\includegraphics[width=1.02\linewidth]{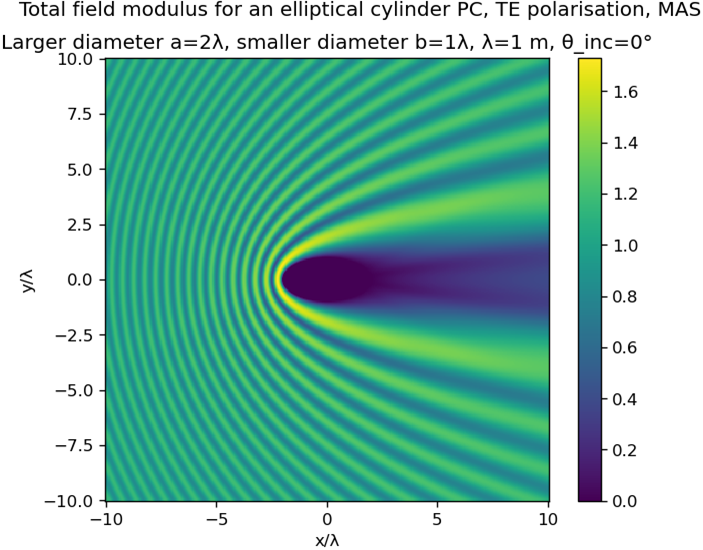}
\caption{Total field modulus for an elliptical cylinder given by the MAS. $q=29.61$, with $m=20$ Mathieu functions for the analytical solution and $n = 40$ sources for the MAS.}
\label{Totalfi}\label{resultfig}
\end{figure}
\begin{figure}[!t]
\includegraphics[scale=0.56]{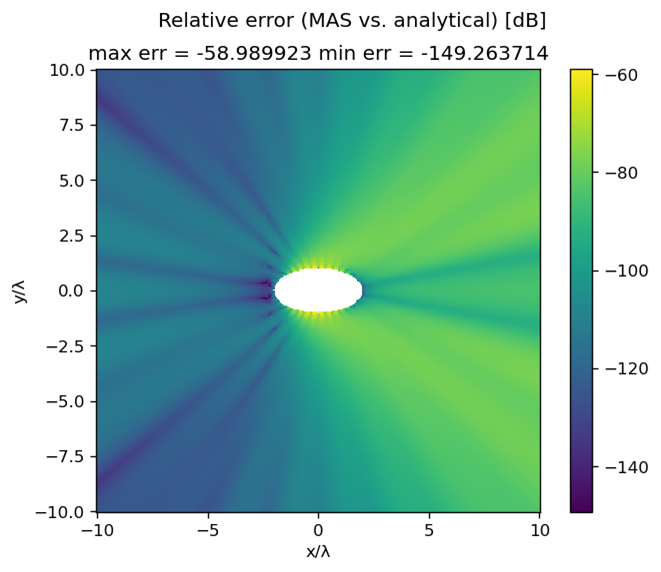}
\caption{Relative error of the MAS compared to the analytical solution. $q=29.61$, with $m=20$ Mathieu functions for the analytical solution and $n = 40$ sources for the MAS.}\label{Relaterr}
\end{figure}
Fig. \ref{fig:ellipse} represents the considered geometry. Sources and collocation points are distributed uniformly, with $\mu_\delta = 0.32\,\mu_s$. \newline
Fig. \ref{Totalfi} shows the total field modulus given by the MAS, including the incident and scattered field, for $n = 40$ auxiliary sources and collocation points. Fig. \ref{Relaterr} displays the relative error between the scattered field from the MAS and the scattered field from the analytical method \cite{van_den_berg_diffraction_1973}. The agreement is excellent, as expected, with a maximum relative error close to $-60$ dB. The analytical solution includes a sum of $20$ even Mathieu functions and $20$ odd Mathieu functions.\newline
The first step of the numerical simulations consists in validating equations \eqref{vpcm} and \eqref{vpsm}. The diagonal coefficients or ``eigenvalues" of the specified rescaled Mathieu compression are considered; a separate full-matrix eigenvalue computation would not be equivalent. To achieve this, the eigenvalues are computed numerically, through the matrix product computation, and compared with the theoretical values. Similarly to the analytical method, $20$ odd eigenvalues and $20$ even eigenvalues are respectively considered for equations \eqref{vpcm} and \eqref{vpsm}. The number of source angles is raised to $40$ to avoid aliasing. Fig. \ref{resultfig} proved this number adequate for the MAS to converge. \newline
Fig. \ref{fig:vpcomp} and \ref{fig:vpcomp2} display this comparison in the even and odd cases. The results show a clear agreement between the numerical values and the theoretical estimates. 


\begin{figure}[!t]
    \centering
    \includegraphics[width=\linewidth]{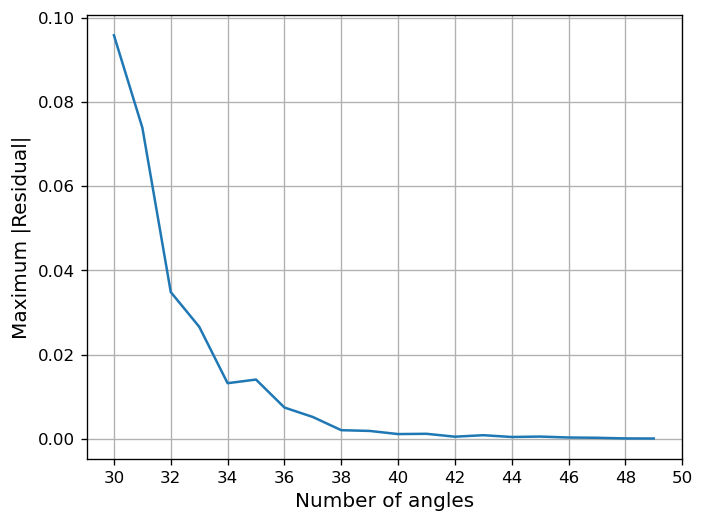}
    \caption{Maximum residual of the impedance matrix, depending on the number of angles. $q =29.61$.}
\label{FG}
\end{figure}

The second step assesses the checkerboard structure and the continuous relation \eqref{resid} through the discrete MAS matrix. Let $\widehat{\mathbf{Z}}^{(P)}=\mat{F}^{\dag}\mat{Z}\mat{F}$. The left-hand side of \eqref{resid} is evaluated after replacing $\widehat{T}_{k,s}$ by $\widehat{\mathbf{Z}}^{(P)}_{k,s}$. Its maximum normalized residual \begin{equation}
   \resizebox{\hsize}{!}{$ \frac{\max_{k,s}\left|(k^2-s^2)\,\hat{\textbf{Z}}_{k,s} + q\bigl( \hat{\textbf{Z}}^{(P)}_{k-2,s} + \hat{\textbf{Z}}^{(P)}_{k+2,s} - \hat{\textbf{Z}}^{(P)}_{k,s-2} - \hat{\textbf{Z}}^{(P)}_{k,s+2} \bigr)\right|}{||\hat{\textbf{Z}}^{(P)}||_2}$}
\end{equation}
is reported in Fig.~\ref{FG}, depending on the number of angles $P$. For a quantitative convergence test, the maximizing index set and denominator norm must be stated and kept consistent as $P$ changes; the figures currently provide a diagnostic rather than a measured convergence rate.
As the number of angular samples increases, the residual decreases, confirming the convergence of the discrete Fourier representation toward the continuous relation \eqref{resid}.\newline
\begin{figure}
\centering
\includegraphics[width=\linewidth]{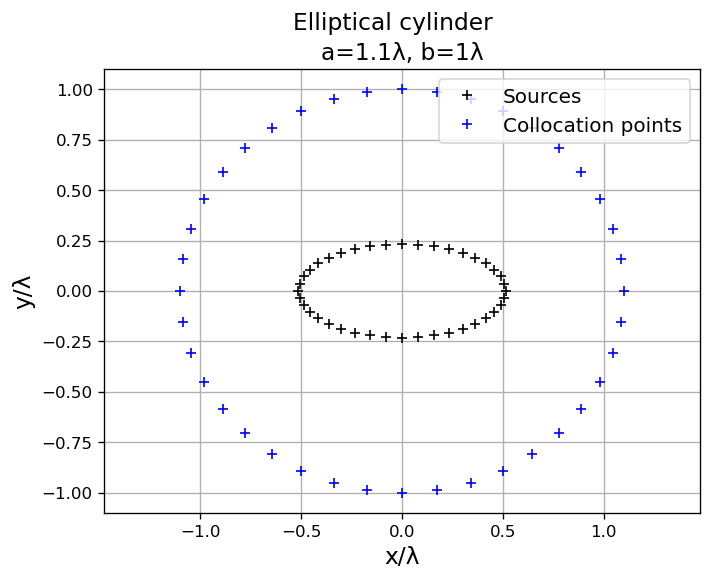} 
\caption{Meshing of the quasi-circular ellipse considered in the numerical simulations.}\label{fig:ellipse2}
\end{figure}
The numerical simulations can also be performed for a quasi-circular ellipse, to assess their evolution. This time, the major semi-axis of the ellipse is $a=1.1\lambda_0$, the minor semi-axis is still $b=\lambda_0$, with $q=2.07$. The confocal ellipse parameter remains $\mu_\delta = 0.32\,\mu_s$, with $n=40$ angle samples and $20$ eigenvalues for $\mat{Z}_1$ and $\mat{Z}_2$. The considered geometry is represented in Fig. \ref{fig:ellipse2}. With this configuration, the relative error is much smaller, with a maximum around $-150$ dB.\newline
\begin{figure}[!t]
\centering
\includegraphics[width=\linewidth]{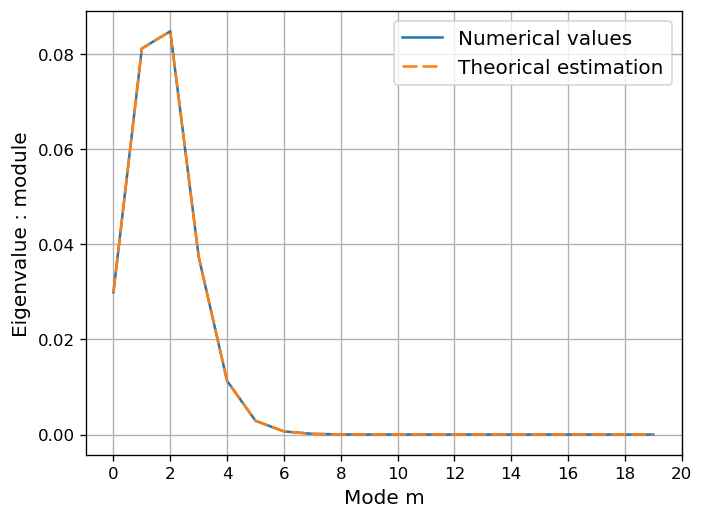}
\caption{Comparison between theoretical and numerical eigenvalues depending on the mode $m$. Even case $\mat{Z}_1$. $\lambda_0=1$ m, $a = 1.1$ m, $b = 1$ m, $q =2.07$.}
\label{fig:vpcompc}
\end{figure}
\begin{figure}
\includegraphics[width=\linewidth]{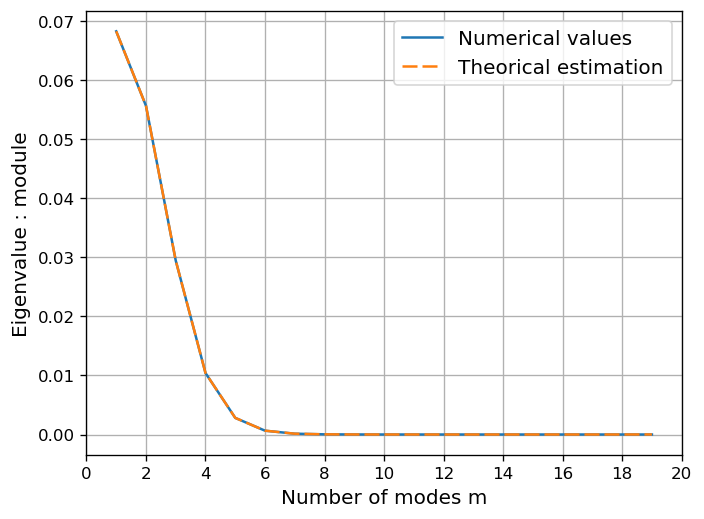}
\caption{Comparison between theoretical and numerical eigenvalues depending on the mode $m$. Odd case $\mat{Z}_2$. $\lambda_0=1$ m, $a = 1.1$ m, $b = 1$ m, $q =2.07$.}
\label{fig:vpcomps}
\end{figure}
The new theoretical and numerical eigenvalues for the reduction of $\mat{Z}_1$ and $\mat{Z}_2$ are shown in Figs. \ref{fig:vpcompc} and \ref{fig:vpcomps}. Both confirm the agreement between the curves. They also display similar curve shapes to the ones of Figs. \ref{fig:vpcomp} and \ref{fig:vpcomp2}, but shifted toward smaller angular values. It means that less eigenvalues are necessary to describe both matrices, exactly as less Mathieu functions would be expected for the analytical method. This is due to the products $\mathrm{Mc}_n^{(1)}(\mu_\delta, q)\,
\mathrm{Mc}_n^{(3)}(\mu_s, q)$ and $\mathrm{Ms}_n^{(1)}(\mu_\delta, q)\,
\mathrm{Ms}_n^{(3)}(\mu_s, q)$, decreasing in both cases. A similar product is used for the analytical solution, thus explaining the analogy.\newline
\begin{figure}[!t]
    \centering
    \includegraphics[width=\linewidth]{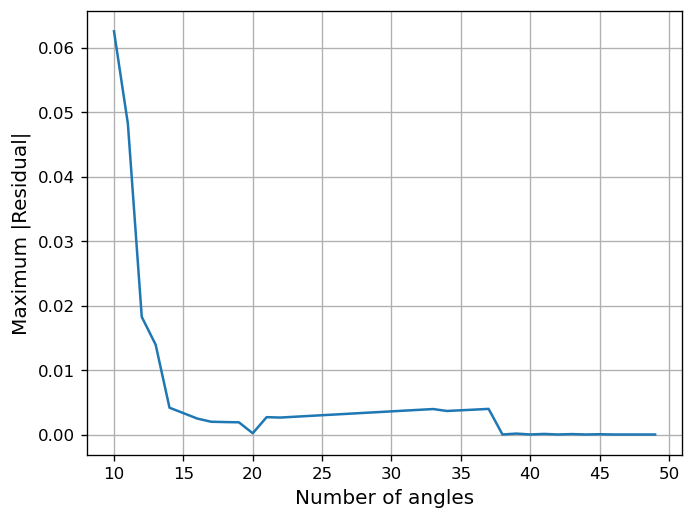}
    \caption{Maximum residual of the impedance matrix, depending on the number of angles. $q =2.07$.}
\label{RC}
\end{figure}
The new normalized and discretized residual of $\hat{\textbf{Z}}^{(P)}$ is represented in Fig. \ref{RC}. This time, the maximum residual is smaller for a much lower number of sources, starting from $10$ auxiliary sources. Although the curve shows some numerical jumps, it exhibits an overall decrease, which confirms the convergence of the residual.\newline
\begin{figure}
\centering
\includegraphics[width=\linewidth]{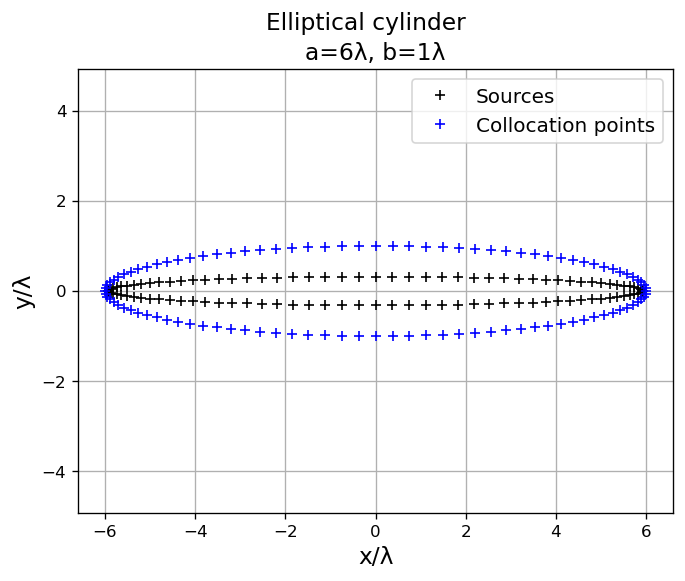} 
\caption{Meshing of the extreme ellipse considered in the numerical simulations.}\label{fig:ellipse3}
\end{figure}
On the other side, we can consider a scatterer with a much more pronounced eccentricity. Fig. \ref{fig:ellipse3} displays an elliptical cylinder with $a=6\lambda_0$, $b=\lambda_0$, and $q=345.44$. The confocal ellipse parameter remains $\mu_\delta = 0.32\,\mu_s$. This time, the number of sources is raised to $n=100$ and $40$ eigenvalues are selected for $\mat{Z}_1$ and $\mat{Z}_2$. As a matter of fact, aliasing would occur for the previous $n=40$ angle values. Due to the extreme configuration, numerical errors arise in the comparison with the analytical method and the maximum relative error is around $-5$ dB.\newline
\begin{figure}[!t]
\centering
\includegraphics[width=\linewidth]{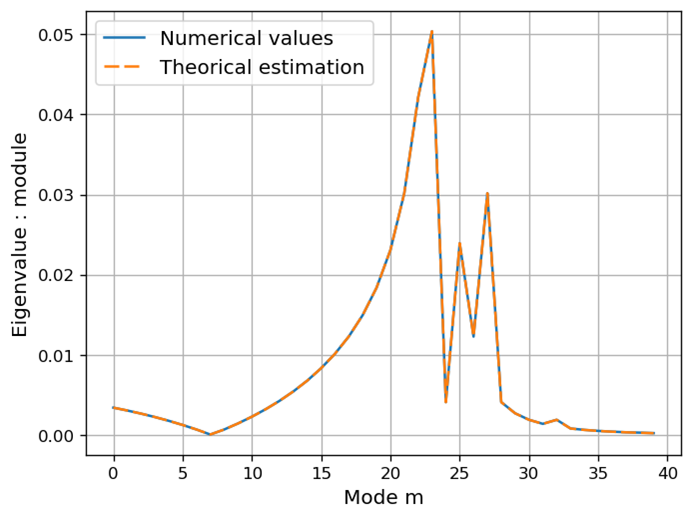}
\caption{Comparison between theoretical and numerical eigenvalues depending on the mode $m$. Even case $\mat{Z}_1$. $\lambda_0=1$ m, $a = 6$ m, $b = 1$ m, $q =345.44$.}\label{fig:vpc3}
\end{figure}
\begin{figure}[!t]
\includegraphics[width=\linewidth]{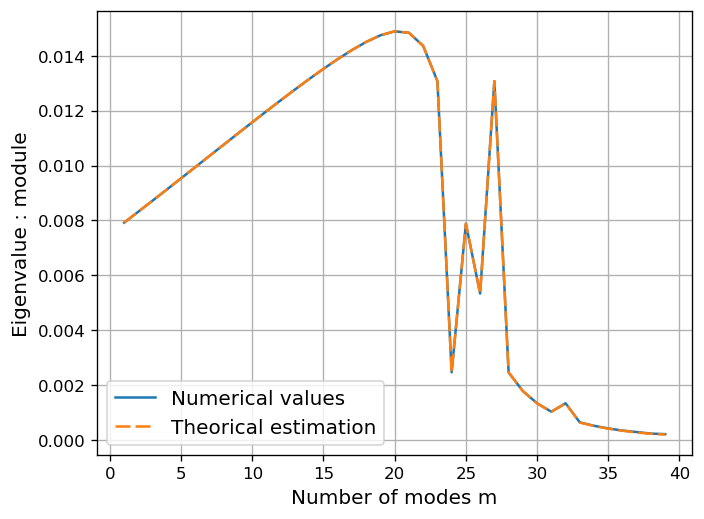}
\caption{Comparison between theoretical and numerical eigenvalues depending on the mode $m$. Odd case $\mat{Z}_2$. $\lambda_0=1$ m, $a = 6$ m, $b = 1$ m, $q =345.44$.}
\label{fig:vps3}
\end{figure}
Figs. \ref{fig:vpc3} and \ref{fig:vps3} represent the comparison between the theoretical and numerical eigenvalues for the even and odd cases. While the numerical results agree with the theoretical values, the eigenvalues oscillate much more than previously (cf. Figs. \ref{fig:vpcomp} and \ref{fig:vpcomp2} for example). The same comparison with $n=90$ sources, not displayed here, showed that the numerical eigenvalues from the higher modes suffer from aliasing, with the emergence of a small difference between the curves. 
\begin{figure}
    \centering
    \includegraphics[width=\linewidth]{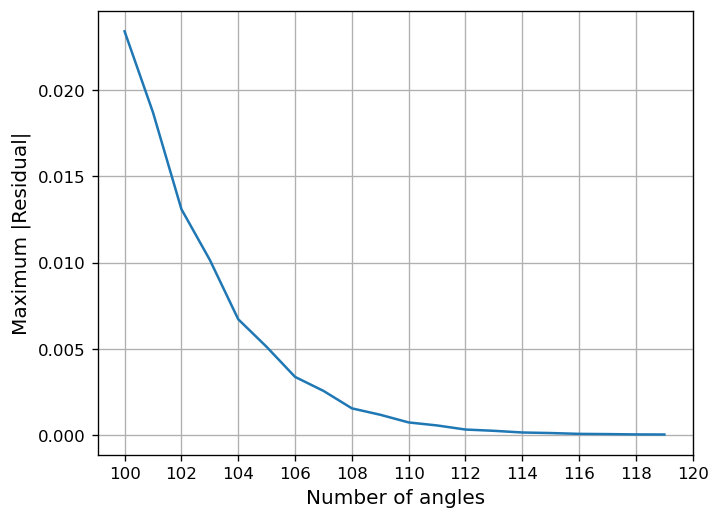}
    \caption{Maximum residual of the impedance matrix, depending on the number of angles. $q =345.44$.}
\label{R3}
\end{figure}
The residual given in Fig. \ref{R3} yields the same conclusion. The maximum residual error for $n=100$ is twice as high as the residual error from Fig. \ref{FG} for $n=40$. These two maxima need not be directly comparable if the set of Fourier indices changes with $P$; a fixed resolved-mode window is preferable for a convergence claim.
Although the residual modulus decreases between $n=100$ and $n=120$, this behavior is an indicator and a good reminder of the susceptibility to aliasing in higher modes.

\section{Conclusion}
In this article, we have calculated the impedance matrix of the MAS for an ellipse, and  given its expansion in Mathieu functions. It allowed us to project it into a reduced diagonal matrix. At the same time, an analysis of the matrix in the Fourier domain revealed a checkerboard structure and relations between every second element. Numerical simulations have been performed for three different ellipses, including a quasi-circular configuration and a high eccentricity geometry. Each evaluation confirmed our theoretical observations. The additional operator estimates establish exact parity on even grids, a discrete aliasing formula, and a controlled band-truncation error in the near-circular regime. The compression and conditioning bounds are mathematical results for specified idealized or truncated operators; assessing their practical benefit across the full parameter range requires dedicated numerical tests.
\section*{Acknowledgments}
This work is supported by the AID (Defence Innovation Agency, Agence de l’Innovation de Défense), a French Defence Ministry entity. The authors would like to express their sincere gratitude.
\section*{Declaration of generative AI and AI-assisted technologies in the manuscript preparation process}
During the preparation of this work the authors used Copilot, Gemini, Deepl for translation purposes, mathematical exposition and Python code assistance. After using this tool/service, the authors carefully edited the content as needed and take full responsibility for the content of the published article.

\bibliographystyle{unsrtnat}
\bibliography{sample}

\end{document}